\documentclass[12pt,letterpaper]{article}

\usepackage{cite}

\begin{document}

\title{Optimization problems constrained by parameter sums}

\author{John C. Nash and Ravi Varadhan}

\date{2024-06-05}

\maketitle

\abstract{%
This article presents a discussion of optimization problems where the objective function \(f(\textbf{x})\) has parameters that are constrained by some scaling, so that \(q(\textbf{x}) = constant\), where this function \(q()\) involves a sum of the parameters, their squares, or similar simple function. Our focus is on ways to use standardized optimization programs to solve such problems rather than specialized codes. Examples are presented with R.
}

\section{Background}\label{background}

We consider problems where we want to minimize or maximize a function
subject to a constraint that the sum of some function of the parameters, e.g.,
their sum or sum of squares, must equal some constant. Because these problems
all have an objective that is dependent on a scaled set of parameters
where the scale is defined by this sum, we will refer to them as \textbf{sumscale}
optimization problems.

We have observed questions about problems like this on the R-help mailing list:

\footnotesize
\begin{verbatim}
Jul 19, 2012 at 10:24 AM, Linh Tran <Tranlm```berkeley.edu> wrote:
> Hi fellow R users,
>
> I am desperately hoping there is an easy way to do this in R.
>
> Say I have three functions:
>
> f(x) = x^2
> f(y) = 2y^2
> f(z) = 3z^2
>
> constrained such that x+y+z=c (let c=1 for simplicity).
>
> I want to find the values of x,y,z that will minimize 
f(x) + f(y) + f(z).
\end{verbatim}
\normalsize

For this particular problem, if the parameters \(x\), \(y\) and \(z\) are non-negative,
this problem can be solved as a Quadratic Program. We revisit this problem at the
end of the article.

Other examples of this type of objective function are:

\begin{itemize}
\item
  The maximum volume of a regular polyhedron where the sum of the lengths
  of the sides is fixed.
\item
  The minimum negative log likelihood for a multinomial model.
\item
  The Rayleigh Quotient for the maximal or minimal eigensolutions of a matrix, where
  the eigenvectors should be normalized so the square norm of the vector is 1.
\item
  The minimum of the extended Rosenbrock function of the form given by
  the \texttt{adagio} package (Borchers (\cite{p-adagio}{2022})) on the unit ball, that is, where
  the sum of squares of the parameters is 1.
\end{itemize}

\section{An illustrative example}\label{an-illustrative-example}

For the moment, let us consider a very basic example, which is

\textbf{Problem A: Minimize \(( - \prod{\textbf{x}})\) subject to \(\sum{\textbf{x}}=1\)}

This is a very simplified version of the multinomial maximum likelihood problem.
We will specify that there are \(nprob\) parameters, though they have only \((nprob - 1)\)
degrees of freedom due to the summation constraint. We have chosen to
specify the problem
as a minimization, but the alternative maximization simply omits the negative sign
in the objective. It is obvious that a solution to this problem has all parameters
at \texttt{1/nprob}. A variant of this problem that gives parameters having unequal
values uses a constraint where the inner product of \(x\) with some vector of
values \(v\) is required to equal some number.

Our goal in this article is to explore how such problems may be most
efficiently approached using general optimization tools rather than
special-purpose software. Moreover, if possible we want to

\begin{itemize}
\item
  develop usable recipes for specifying the optimization to solve such
  problems for general users rather than mathematical specialists
\item
  learn what types of optimization solvers are most effective for
  different cases, or which are to be avoided
\item
  discover whether it is useful to apply bounds constraints on parameters.
\end{itemize}

\section{Difficulties using general optimization with sumscale problems}\label{difficulties-using-general-optimization-with-sumscale-problems}

The basic example illustrates some of the main issues with sumscale problems.

First, we may encounter \textbf{scaling} difficulties. Noting that the solution of
our problem is at \(1/nprob\), the objective function value at the solution is

\[   -  nprob ^ {(- nprob)} \]

The smallest number R can represent is approximately \texttt{2e-308} of which the
base 2 logarithm
is approximately \texttt{-708}, which the objective function solution achieves at
approximately \(nprob = 142\). This is not an unusually large number of
parameters in practice.

The second issue, and the one that causes more difficulties for software,
is that any solution is unique only to within an arbitrary
positive scaling factor. This complicates optimization, since convergence
or termination tests check if the parameters have changed. Clearly, we can change
the parameter values of the original objective function by a constant factor
without changing the objective function value.

\subsection{Preliminary attempted solution}\label{preliminary-attempted-solution}

Let us use the basic example above to consider how we might formulate it
for a computational solution.

Surprisingly, occasionally we may be able to get away with simply
putting the constraint into
the objective by scaling the parameters as in the following example.

\tiny 
\begin{verbatim}
library(optimx)
pr0 <- function(x) {
- prod(x/sum(x))
}
\end{verbatim}

\begin{verbatim}
#> test the simple product for n=5
\end{verbatim}

\begin{verbatim}
#>                    p1 s1        p2 s2        p3 s3        p4 s4        p5 s5
#> nvm         0.1459940    0.1459940    0.1459940    0.1459940    0.1459940   
#> ncg         0.1406366    0.1406421    0.1406365    0.1406408    0.1406431   
#> Nelder-Mead 0.1311841    0.1311930    0.1311787    0.1311864    0.1311869   
#> BFGS        0.1462928    0.1463680    0.1462937    0.1463611    0.1462971   
#>                value fevals gevals hevals conv kkt1  kkt2 xtime
#> nvm         -0.00032     40     33      0    0 TRUE FALSE 0.003
#> ncg         -0.00032      9      7      0    0 TRUE FALSE 0.001
#> Nelder-Mead -0.00032    228      0      0    0 TRUE FALSE 0.003
#> BFGS        -0.00032     16     14      0    0 TRUE FALSE 0.001
\end{verbatim}

\begin{verbatim}
#> Method:  Nelder-Mead  max(abs(error_in_parameters)) = 1.094366e-05 
#> Method:  BFGS  max(abs(error_in_parameters)) = 6.215476e-05 
#> Method:  ncg  max(abs(error_in_parameters)) = 4.779663e-06 
#> Method:  nvm  max(abs(error_in_parameters)) = 1.017677e-10
\end{verbatim}
\normalsize

This approach does surprisingly well, but we do note considerable variation in the
sizes of the (raw) parameters, though it is easy to scale them.

Another possibility is to select one of the parameters and solve for it in
terms of the others. Let this
be the last parameter \(x_n\), so that the set of parameters to be
optimized is \(\textbf{y} = (x_1, x_1, ..., x_{n-1})\) where
\(n\) is the original size of our problem. We now have the unconstrained problem

\[ minimize ( - (\prod{\textbf{y}}) * (1 - \sum{y} ) ) \]

This is easily coded and tried. We will use a very simple start, namely, the sequence
\(1,2, ..., (n-1)\) scaled by \(1/n^2\). We will also specify that the gradient is to be
computed by a central approximation (see \texttt{gcentral.R} from package \textbf{optimx}).

\tiny
\begin{verbatim}
pr1 <- function(y) {
- prod(y)*(1-sum(y))
}
\end{verbatim}

\begin{verbatim}
#> test the simple product for n=5
\end{verbatim}

\begin{verbatim}
#>                    p1 s1        p2 s2        p3 s3        p4 s4    value fevals
#> nvm         0.2000000    0.2000000    0.2000000    0.2000000    -0.00032    124
#> BFGS        0.2000000    0.1999986    0.2000037    0.1999985    -0.00032    102
#> Nelder-Mead 0.2000034    0.1999983    0.2000017    0.2000021    -0.00032    331
#> ncg         0.1999934    0.1999958    0.2000014    0.2000129    -0.00032     20
#>             gevals hevals conv kkt1 kkt2 xtime
#> nvm            109      0    3 TRUE TRUE 0.006
#> BFGS            96      0    0 TRUE TRUE 0.003
#> Nelder-Mead      0      0    0 TRUE TRUE 0.002
#> ncg             10      0    0 TRUE TRUE 0.000
\end{verbatim}

\begin{verbatim}
#> Method:  Nelder-Mead  max(abs(error_in_parameters)) = 5.544115e-06 
#> Method:  BFGS  max(abs(error_in_parameters)) = 3.715673e-06 
#> Method:  ncg  max(abs(error_in_parameters)) = 1.289054e-05 
#> Method:  nvm  max(abs(error_in_parameters)) = 2.24181e-09
\end{verbatim}
\normalsize

While these codes work well for small \(n\), it is fairly easy to see that there are
computational problems as the size of the problem increases. We have already noted
that the function will underflow when \(n\approx142\).

\subsubsection{Another naive formulation}\label{another-naive-formulation}

One of the features added in some of the solvers for \texttt{optimx} is the possibility
of \textbf{masks}, that is, fixed parameters. We want to fix just one parameter, but
doing so fixes the scaling. This will alter the product value unless we scale
to the summation constraint, which we do in function \texttt{pr0} below. Unfortunately,
gradient methods seem to do very poorly on this function, even for small \(n\).

Indeed, in most cases, the gradient methods made no progress towards a solution.

\subsection{Avoiding the poor function scaling}\label{avoiding-the-poor-function-scaling}

Traditionally, statisticians solve maximum likelihood problems by \textbf{minimizing}
the negative log-likelihood. That is, the objective function is formed as (-1) times
the logarithm of the likelihood. This converts our product to a sum. Moreover, the gradient
can be developed more easily. Let us choose to eliminate the last parameter using
the summation constraint. Note that an original problem of size \texttt{nprob}
has \texttt{(nprob\ -\ 1)} working parameters. As programs that try to find the minimum may change
the parameters so that logarithms of non-positive numbers are attempted, we have put
some safeguards in the function \texttt{nll()}. At this point we have assumed the gradient
calculation is only attempted if the function can be computed satisfactorily, so we
have not put safeguards in the gradient. We could also put bounds on the parameters,
since they cannot be non-positive, nor exceed 1. This is, of course, extra work
for the user.

\tiny
\begin{verbatim}
# nll uses (nprob - 1) parameters,  y <- x[1:(nprob - 1)]
# implicit last param is (1 - sum(rest)); nll is fn of y, not x
nll <- function(y) { # safeguarded "solve for 1" approach
  if ((any(y <= 10*.Machine$double.xmin)) || (sum(y)>1-.Machine$double.eps))
    .Machine$double.xmax
  else   - sum(log(y)) - log(1-sum(y))
}
nll.g <- function(y) { - 1/y + 1/(1-sum(y))} # not safeguarded
\end{verbatim}
\normalsize

Let us now attempt solutions for problems in 5, 100, and 1000
parameters with a selection
of solvers and different constraints. Note that the solvers ``optim::L-BFGS-B'' and
``lbfgsb3c'' both failed to solve these problems, and we do not report results.
Since we have an analytical gradient, we will use it. However, use of

\tiny
\begin{verbatim}
gr="grcentral"
\end{verbatim}
\normalsize

\noindent
gave quite similar (and acceptable) results using the central difference approximation.
The use of bounds (preferably not too tight) seems from a few examples to improve
performance by keeping parameters in a reasonable domain.

\tiny
\begin{verbatim}
#> NOTE: using a forward difference gave check errors on the gradient function
\end{verbatim}

\begin{verbatim}
#> n=5; analytical gradient with [0,1] bounds
\end{verbatim}

\begin{verbatim}
#>               p3 s3        p4 s4   value fevals gevals hevals conv kkt1 kkt2 xtime
#> spg    0.2000000    0.2000000    8.04719     52     15      0    0   NA   NA 0.005
#> nlminb 0.1999989    0.1999992    8.04719     24     12      0    0   NA   NA 0.002
#> nvm    0.2000000    0.2000000    8.04719     28     14      0    0   NA   NA 0.001
#> ncg    0.2000000    0.2000000    8.04719     26     12      0    0   NA   NA 0.001
#> tnewt  0.2000000    0.2000000    8.04719     20     19      0    0   NA   NA 0.002
\end{verbatim}

\begin{verbatim}
#> Method:  spg  max(abs(error_in_parameters)) = 4.982131e-08   fval= 8.04719 
#> Method:  nlminb  max(abs(error_in_parameters)) = 2.580296e-06   fval= 8.04719 
#> Method:  nvm  max(abs(error_in_parameters)) = 7.311371e-10   fval= 8.04719 
#> Method:  ncg  max(abs(error_in_parameters)) = 6.850362e-09   fval= 8.04719 
#> Method:  tnewt  max(abs(error_in_parameters)) = 1.311518e-10   fval= 8.04719
\end{verbatim}

\begin{verbatim}
#> n=100, analytical gradient, [0,1] bounds
\end{verbatim}

\begin{verbatim}
#>               p98 s98        p99 s99   value fevals gevals hevals conv kkt1 kkt2 xtime
#> spg    0.01000000     0.01000000     460.517    890     75      0    0   NA   NA 0.016
#> nlminb 0.01001132     0.01001328     460.517    170    151      0    1   NA   NA 0.008
#> nvm    0.01000000     0.01000000     460.517    400    123      0    0   NA   NA 0.048
#> ncg    0.01000000     0.01000000     460.517    114     43      0    0   NA   NA 0.013
#> tnewt  0.01000000     0.01000000     460.517     55     54      0    0   NA   NA 0.003
\end{verbatim}

\begin{verbatim}
#> Method:  spg  max(abs(error_in_parameters)) = 4.733013e-11   fval= 460.517 
#> Method:  nlminb  max(abs(error_in_parameters)) = 1.327973e-05   fval= 460.517 
#> Method:  nvm  max(abs(error_in_parameters)) = 9.068386e-09   fval= 460.517 
#> Method:  ncg  max(abs(error_in_parameters)) = 2.155583e-10   fval= 460.517 
#> Method:  tnewt  max(abs(error_in_parameters)) = 4.53821e-14   fval= 460.517
\end{verbatim}

\begin{verbatim}
#> n=1000, analytical gradient, [0,1] bounds, maxit=1000
\end{verbatim}

\begin{verbatim}
#> Error in if (!max.diff < checkGrad.tol) { : 
#>   missing value where TRUE/FALSE needed
\end{verbatim}

\begin{verbatim}
#>               p998 s998        p999 s999         value fevals gevals hevals conv kkt1
#> spg             NA               NA      8.988466e+307   7995      1      0 9999   NA
#> nlminb 0.001005072      0.001005056       6.913893e+03    187    151      0    1   NA
#> nvm    0.001000000      0.001000000       6.907755e+03   1766    515      0    0   NA
#> ncg    0.001000000      0.001000000       6.907755e+03    312     99      0    0   NA
#> tnewt  0.001000000      0.001000000       6.907755e+03    148    147      0    0   NA
#>        kkt2 xtime
#> spg      NA 0.219
#> nlminb   NA 0.576
#> nvm      NA 9.904
#> ncg      NA 0.334
#> tnewt    NA 0.010
\end{verbatim}

\begin{verbatim}
#> Method  spg  failed  -- fval= 8.988466e+307 
#> Method:  nlminb  max(abs(error_in_parameters)) = 0.007087251   fval= 6913.893 
#> Method:  nvm  max(abs(error_in_parameters)) = 1.514042e-09   fval= 6907.755 
#> Method:  ncg  max(abs(error_in_parameters)) = 7.513133e-09   fval= 6907.755 
#> Method:  tnewt  max(abs(error_in_parameters)) = 1.138514e-11   fval= 6907.755
\end{verbatim}
\normalsize

We note quite different performance for different methods. For example, \texttt{nvm} has
a much longer timing than the rest. All methods here have done
reasonably well for accuracy of the solution, though \texttt{spg} has not produced
a solution for \texttt{n\ =\ 1000} because a programmatic safety
check on the analytic gradient has falsely suggested there is an error. (Method
\texttt{spg} was never applied in this case.)

\subsubsection{Parameter scaling}\label{parameter-scaling}

In contrast to the explicit use of the summation to get the last parameter, let
us instead try minimizing a function of scaled parameters using methods
\texttt{nvm} and \texttt{ncg}. In some trials we mask one of the parameters.

\tiny
\begin{verbatim}
#> Function:  (-1)*sum(log(x/sum(x)))
\end{verbatim}

\begin{verbatim}
#> n=5:
\end{verbatim}

\begin{verbatim}
#>            p3 s3        p4 s4        p5 s5   value fevals gevals hevals conv kkt1 kkt2
#> ncg 0.4122608    0.4122608    0.4122608    8.04719     11      7      0    0   NA   NA
#> nvm 0.8691075    0.8691075    0.8691075    8.04719     20     12      0    0   NA   NA
#>     xtime
#> ncg 0.000
#> nvm 0.002
\end{verbatim}

\begin{verbatim}
#> Method:  ncg  max(abs(error_in_parameters)) = 1.932914e-08   fval= 8.04719 
#> Method:  nvm  max(abs(error_in_parameters)) = 1.598058e-10   fval= 8.04719
\end{verbatim}

\begin{verbatim}
#> n=100:
\end{verbatim}

\begin{verbatim}
#>           p98 s98       p99 s99      p100 s100   value fevals gevals hevals conv kkt1
#> ncg 0.6188335     0.6188335     0.6188334      460.517     13     11      0    0   NA
#> nvm 0.6297044     0.6297044     0.6297044      460.517     35     21      0    0   NA
#>     kkt2 xtime
#> ncg   NA 0.003
#> nvm   NA 0.007
\end{verbatim}

\begin{verbatim}
#> Method:  ncg  max(abs(error_in_parameters)) = 1.436611e-08   fval= 460.517 
#> Method:  nvm  max(abs(error_in_parameters)) = 1.004048e-09   fval= 460.517
\end{verbatim}

\begin{verbatim}
#> n=100 masked last parameter:
\end{verbatim}

\begin{verbatim}
#>            p98 s98        p99 s99 p100 s100   value fevals gevals hevals conv kkt1 kkt2
#> ncg 0.02000001     0.02000001     0.02    M 460.517      9      4      0    0   NA   NA
#> nvm 0.02000000     0.02000000     0.02    M 460.517     45      9      0    0   NA   NA
#>     xtime
#> ncg 0.002
#> nvm 0.006
\end{verbatim}

\begin{verbatim}
#> Method:  ncg  max(abs(error_in_parameters)) = 4.461748e-09   fval= 460.517 
#> Method:  nvm  max(abs(error_in_parameters)) = 4.858782e-10   fval= 460.517
\end{verbatim}

\begin{verbatim}
#> n=1000 masked:
\end{verbatim}

\begin{verbatim}
#>           p998 s998       p999 s999 p1000 s1000    value fevals gevals hevals conv kkt1
#> ncg 0.01999512      0.01999512       0.02     M 6907.755      6      3      0    0   NA
#> nvm 0.02000000      0.02000000       0.02     M 6907.755     26      7      0    0   NA
#>     kkt2 xtime
#> ncg   NA 0.007
#> nvm   NA 0.102
\end{verbatim}

\begin{verbatim}
#> Method:  ncg  max(abs(error_in_parameters)) = 2.43707e-07   fval= 6907.755 
#> Method:  nvm  max(abs(error_in_parameters)) = 2.306793e-11   fval= 6907.755
\end{verbatim}
\normalsize

When fixing a parameter with a mask, we must scale the solution parameters before
comparing them to the known solution. With or without masking, both methods do quite
well. However, for the larger problem, the storage of the
approximate inverse Hessian in ``nvm'' occasions much greater computational
effort.

We can easily try several other optimization methods using the \texttt{optimx} package.
The calls are given in the example below for \texttt{n\ =\ 100} on the \texttt{nll()} function
using analytical gradient \texttt{nll.g()}. Note that we do not ask for
\texttt{method="ALL"} or \texttt{method="MOST"} as we found that some of the methods,
in particular those using Powell's
quadratic approximation methods, seem to get ``stuck''. This may be a
consequence of some bugs recently discovered in these methods, see Ragonneau and Zhang (\cite{PDFO23}{2023}).
We selected the methods
``spg'', ``nlminb'', ``nvm'', ``ncg'', and ``tnewt'' which use a variety of algorithms.
(``nvm'' and ``ncg'' are repeated for ease of comparison.)
All handle bounds constraints, which in this example seems to confer an advantage
both for solution accuracy and performance.

\tiny
\begin{verbatim}
#> no bounds
\end{verbatim}

\begin{verbatim}
#>                p98 s98         p99 s99    value fevals gevals hevals conv kkt1 kkt2 xtime
#> spg    0.009999999     0.009999999     460.5170    969    152      0    0   NA   NA 0.016
#> nlminb 0.010000000     0.010000000     460.5170    178    128      0    0   NA   NA 0.007
#> nvm    0.010000000     0.010000000     460.5170    742    137      0    0   NA   NA 0.026
#> ncg    0.010000000     0.010000000     460.5170    111     41      0    0   NA   NA 0.001
#> tnewt  0.009900000     0.010000000     548.7873     14     13      0 9999   NA   NA 0.001
\end{verbatim}

\begin{verbatim}
#> Method:  spg  max(abs(error_in_parameters)) = 2.590979e-07   fval= 460.517 
#> Method:  nlminb  max(abs(error_in_parameters)) = 3.531403e-08   fval= 460.517 
#> Method:  nvm  max(abs(error_in_parameters)) = 1.082334e-09   fval= 460.517 
#> Method:  ncg  max(abs(error_in_parameters)) = 2.501705e-10   fval= 460.517 
#> Method  tnewt  failed  -- fval= 548.7873
\end{verbatim}

\begin{verbatim}
#> simple [0,1] bounds
\end{verbatim}

\begin{verbatim}
#>               p98 s98        p99 s99   value fevals gevals hevals conv kkt1 kkt2 xtime
#> spg    0.01000000     0.01000000     460.517    890     75      0    0   NA   NA 0.015
#> nlminb 0.01001132     0.01001328     460.517    170    151      0    1   NA   NA 0.008
#> nvm    0.01000000     0.01000000     460.517    400    123      0    0   NA   NA 0.051
#> ncg    0.01000000     0.01000000     460.517    114     43      0    0   NA   NA 0.014
#> tnewt  0.01000000     0.01000000     460.517     55     54      0    0   NA   NA 0.002
\end{verbatim}

\begin{verbatim}
#> Method:  spg  max(abs(error_in_parameters)) = 4.733013e-11   fval= 460.517 
#> Method:  nlminb  max(abs(error_in_parameters)) = 1.327973e-05   fval= 460.517 
#> Method:  nvm  max(abs(error_in_parameters)) = 9.068386e-09   fval= 460.517 
#> Method:  ncg  max(abs(error_in_parameters)) = 2.155583e-10   fval= 460.517 
#> Method:  tnewt  max(abs(error_in_parameters)) = 4.53821e-14   fval= 460.517
\end{verbatim}
\normalsize

Oddly, method ``nlminb'' gets a poorer result when bounds are imposed.

We tried the same examples with numerically approximated gradients. For simple
forward differences, we did not get satisfactory solutions. Central differences
gave reasonable solutions, but generally the parameters had slightly greater
deviations from the known solution than with analytic gradients. We note that
numerical gradient routines evaluate the objective at slightly modified parameter
values in order to obtain the approximations, but as far as we are aware no
such routines check that bounds are respected. Such checks would
almost certainly degrade
performance, which is already a concern when using numerical gradient approximations.
Note also that a lower bound of 0 on parameters is not adequate,
since \(log(0)\) is undefined. Choosing a
bound large enough to avoid the logarithm of a zero or negative argument
while still being small enough to allow for parameter optimization is a non-trivial matter.

A common recommendation for solving bounded function minimization problems with R is
the solver ``L-BFGS-B'' from base-R \texttt{optim()} (Byrd et al. (\cite{Byrd95}{1995})). However,
it fails to find a solution for function \texttt{nll()},
in that it tries to compute this function with parameters that are inadmissible.
We have not determined the exact cause of this failure, but observe that
the particular form of \(log(1-sum(x))\)
is undefined once the argument of the logarithm is negative. Indeed, this
is the basis of logarithmic barrier functions for constraints. There
is a similar issue with any of the \(nprob - 1\) parameters near to or less than
zero.

\section{Transformed problems or parameters}\label{transformed-problems-or-parameters}

When problems give difficulties, it is common to re-formulate them by
transformations of the function or the parameters.

\subsection{log() transformation of parameters}\label{log-transformation-of-parameters}

A common method to ensure parameters are positive is to transform
them. In the present case, optimizing over
parameters that are the logarithms of the parameters above
ensures we have positive arguments to most of the
elements of the negative log likelihood. The code in
the following example optimizes the log parameters
\texttt{lx} and not \texttt{x}, so we must back-transform to get the
solution in terms of the original parameters.

\tiny
\begin{verbatim}
enll <- function(lx) {
    x<-exp(lx)
    fval<-  - sum( log( x/sum(x) ) ) 
}
enll.g <- function(lx){
    x<-exp(lx)
    g<-length(x)/sum(x) - 1/x
    gval<-g*exp(lx)
}
\end{verbatim}
\normalsize

As with our function \texttt{pr0()}, the constraint is subsumed into the scaling by \texttt{sum(x}).

\tiny
\begin{verbatim}
#> Initial objective function = 464.6491
\end{verbatim}

\begin{verbatim}
#>          p98 s98   p99 s99  p100 s100   value fevals gevals hevals conv kkt1  kkt2 xtime
#> spg    0.505     0.505     0.505      460.517    808      7      0    0 TRUE FALSE 0.009
#> nlminb 0.505     0.505     0.505      460.517      8      6      0    0 TRUE FALSE 0.001
#> nvm    0.505     0.505     0.505      460.517     16      6      0    0 TRUE FALSE 0.001
#> ncg    0.505     0.505     0.505      460.517      9      6      0    0 TRUE FALSE 0.000
#> tnewt  0.505     0.505     0.505      460.517     13     12      0    0 TRUE FALSE 0.001
\end{verbatim}

\begin{verbatim}
#> Method:  spg  max(abs(error_in_parameters)) = 3.578494e-10 
#> Method:  nlminb  max(abs(error_in_parameters)) = 9.061594e-10 
#> Method:  nvm  max(abs(error_in_parameters)) = 2.352231e-11 
#> Method:  ncg  max(abs(error_in_parameters)) = 3.74142e-12 
#> Method:  tnewt  max(abs(error_in_parameters)) = 8.276366e-15
\end{verbatim}
\normalsize

While there are reasons to think that the indeterminacy of the parameters
might upset the optimization codes in this approach, it seems that in practice,
the objective and gradient here are generally well-behaved.
We do, however, note that the objective function at the start and end of the
optimization differs by only a small amount.

\subsection{`Leave one out' approach with log parameters}\label{leave-one-out-approach-with-log-parameters}

Leaving out the last parameter and computing it from the other
\texttt{nprob\ -\ 1} values is also possible. The parameters are now in
log form, which we must remember when submitting starting values.
With log form parameters, we can also apply bounds, but need to note
that the values will be such that the unlogged parameters make
sense. This approach also complicates the gradient slightly.

\tiny
\begin{verbatim}
xnll <- function(lx) { # last parameter determined from rest
  x<-exp(lx)
  xn <- 1.0 - sum(x)
  if (xn < 1e-8) {
    fval <- 1e+25 # invalid parameters -- can't have negative final parameter
  } else {
#    xplus <- c(x, xn) 
    fval <- - sum(log(x)) - log(xn)
  }
}
xnll.g <- function(lx) {
  x <- exp(lx)
  xn <- 1.0 - sum(x)
  gx <- -1/x + 1/xn
  gg <- gx*x
}
\end{verbatim}

\begin{verbatim}
#> Expected optimal value = 460.517
\end{verbatim}

\begin{verbatim}
#> Initial enll function value = 1089.002
\end{verbatim}

\begin{verbatim}
#> Initial xnll function value = 1089.002
\end{verbatim}

\begin{verbatim}
#>             p98 s98      p99 s99   value fevals gevals hevals conv kkt1 kkt2 xtime
#> nlminb -4.60517     -4.60517     460.517     34     19      0    0 TRUE TRUE 0.005
#> nvm    -4.60517     -4.60517     460.517     37     16      0    0 TRUE TRUE 0.006
#> tnewt  -4.60517     -4.60517     460.517     25     24      0    0 TRUE TRUE 0.002
#> ncg    -4.60517     -4.60517     460.517     32     24      0    0 TRUE TRUE 0.005
#> spg    -4.60517     -4.60517     460.517    820     16      0    0 TRUE TRUE 0.011
\end{verbatim}

\begin{verbatim}
#> Method:  spg  max(abs(error_in_parameters)) = 1.266977e-08 
#> Method:  nlminb  max(abs(error_in_parameters)) = 3.700956e-10 
#> Method:  nvm  max(abs(error_in_parameters)) = 1.088263e-09 
#> Method:  ncg  max(abs(error_in_parameters)) = 3.71333e-09 
#> Method:  tnewt  max(abs(error_in_parameters)) = 6.213466e-12
\end{verbatim}
\normalsize

We note that the unlogged parameter values obtained by optimizing \texttt{xnll()}
rather than \texttt{enll()} are not quite as close to the known solution.
On the other hand, the Kuhn-Karush-Tucker curvature test succeeds
for optimizing \texttt{xnll()}.

\subsection{Using a projection}\label{using-a-projection}

The idea behind projection is to find a point in the domain of feasibility
that is nearest to the given point. Suppose \(x'\) is the candidate value
for the unconstrained optimization, we find its projection \(x\) such that
it satisfies the constraint \(g(x) = 0,\) and simultaneously minimizes
\(||x - x'||.\) The projection feature in \texttt{spg} is useful for
constraints which allow simple projection. Examples of such constraints
are \textit{unit-simplex}, \(\sum x_i = 1 \, (0 \le x_i \le 1)\), and
\textit{unit-length}, \(\sum x^2 = 1\), constraints. Dykstra's algorithm
is a celebrated example of a projection algorithm for finding the
intersection of a set of linear inequalities. \texttt{spg} allows
users to solve linear-inequality constrained optimization problems
via projection.

Objective functions defined by
\((-1)*\prod{\textbf{x}}\) or \((-1)*\sum{log(\textbf{x})}\) will change
with the scale of the parameters. However, the constraint \(\sum{\textbf{x}}=1\)
effectively imposes the scaling
\[ \textbf{x}_{scaled} = \textbf{x}/\sum{\textbf{x}}\]
The optimizer \texttt{spg} from package \texttt{BB} allows us to project our search
direction to satisfy constraints. It allows a gradient function to be supplied,
or an internal gradient approximation can be used.

\tiny
\begin{verbatim}
require(BB, quietly=TRUE)
nllrv <- function(x) {- sum(log(x))}
nllrv.g <- function(x) {- 1/x }
proj <- function(x) {x/sum(x)}
n <- 100
aspg <- spg(par=(1:n)/n^2, fn=nllrv, gr=nllrv.g, project=proj, quiet=TRUE)

aspgn <- spg(par=(1:n)/n^2, fn=nllrv, project=proj, quiet=TRUE)
# using internal grad approx.
\end{verbatim}

\begin{verbatim}
#> F_optimal: with gradient= 460.517   num. grad. approx.= 460.517
\end{verbatim}

\begin{verbatim}
#> fbest =  460.517   when all parameters =  0.01
\end{verbatim}

\begin{verbatim}
#> deviations:  with gradient= 5.055715e-07    num. approx.= 5.055715e-07
\end{verbatim}

\begin{verbatim}
#> Max(abs(deviation from best parameters)):
#>   analytic grad = 5.055715e-07     internal approx. grad= 7.006446e-08
\end{verbatim}
\normalsize

Here the projection \texttt{proj} is the key to success of method
\texttt{spg}. Why \texttt{spg()} does (for this case) better using the internal
numerical gradient approximation than the analytic gradient is an
open question.

Other methods (as yet) do not have the flexibility to impose
the projection directly. We would need to carefully build the projection into
the function(s) and/or the method codes.
This was done by Geradin (\cite{Geradin71}{1971}) for the Rayleigh quotient
problem, but requires a number of changes to the program code.

\subsection{Use of the gradient equations}\label{use-of-the-gradient-equations}

Another approach is to ``solve'' the gradient equations. We can do this with
a sum of squares minimizer, though the \texttt{nls} function in R is
specifically NOT useful as it cannot, by default, deal
with small or zero residuals. However, \texttt{nlfb}
from package \texttt{nlsr} is capable of dealing
with such problems. Unfortunately, it will be slow as it has to
generate the Jacobian by numerical
approximation unless we can provide a function to prepare the
Jacobian analytically. Moreover,
the determination of the Jacobian is still subject to
the unfortunate scaling issues we have
been confronting throughout this article.

We have yet to try this approach.

\section{The Rayleigh Quotient}\label{the-rayleigh-quotient}

The maximal and minimal eigensolutions of a symmetric matrix \(A\) are extrema of
the Rayleigh Quotient

\[ R(x) =  (x' A x)  / (x' x) \]

We can also deal with generalized eigenproblems of the form

\[A x = e B x\]

where B is symmetric and positive definite by using the Rayleigh Quotient

\[ R_g(x) =  (x' A x)  / (x' B x) \]

Once again, the objective is scaled by a function of the parameters, this time
by their sum of squares. Alternatively,
we may think of requiring the \textbf{normalized} eigensolution, which is given as
\[ x_{normalized} = x/sqrt(x' x) \]
We will first try the projected gradient method \texttt{spg} from \texttt{BB}.
Below is the code, where our test uses
a matrix called the Moler matrix Nash (\cite{cnm79}{1979}, Appendix 1). This matrix is simple
to generate and is positive definite, but has one small eigenvalue that may
not be computed to high relative precision. That is, we may get a
number that is precise relative to the largest eigenvalue, but having few
trustworthy digits.

Let us set up some infrastructure for the Rayleigh Quotient of the
matrix. The \texttt{molerbuild(n)} function creates a Moler matrix of
order \texttt{n}. \texttt{raynum(x,\ A)} computes the \textbf{unscaled} quadratic form
computed with \texttt{as.numeric((t(x)\%*\%A)\%*\%x)}, while \texttt{RQ(x,A)} scales
this quantity by dividing it by \texttt{sum(x*x)}.

In order to use the projected gradient method \texttt{BB::spg}, we need a projection
function. The squared nature of the sumscale constraint leads to an indeterminacy
of sign in the solution parameters. That is, if there is a solution \texttt{xstar}, then
\texttt{(-1)*xstar} is also a solution. To get round this, we apply the sign of the first
element of the vector of parameters in the projection.

\tiny
\begin{verbatim}
molerbuild<-function(n){ # Create the moler matrix of order n
   # A[i,j] = i for i=j, min(i,j)-2 otherwise
   A <- matrix(0, nrow = n, ncol = n)
   j <- 1:n
   for (i in 1:n) {
      A[i, 1:i] <- pmin(i, 1:i) - 2
   }
   A <- A + t(A)
   diag(A) <- 1:n
   A
}

raynum<-function(x, A){
   raynum<-as.numeric((t(x)%*%A)%*%x)
}

RQ <- function(x, A){
  RQ <- raynum(x, A)/sum(x*x)
}

# proj<-function(x) { x/sqrt(sum(x*x)) } # Original
proj <- function(x) {sign(x[1]) * x/sqrt(c(crossprod(x)))}
\end{verbatim}
\normalsize

We can compute all the eigensolutions of the moler matrix of order
10 using the R function \texttt{eigen()} which is highly efficient.
However, let us try to use \texttt{BB::spg()} to find the largest eigenvalue by
minimizing the Rayleigh Quotient of \(-A\).

\tiny
\begin{verbatim}
#> Compute eigensolutions of n=10 moler matrix with eigen()
\end{verbatim}

\begin{verbatim}
#> Time for running eigen =  0.2242768  ms
\end{verbatim}

\begin{verbatim}
#> Maximal eigenvalue = 31.58981  RQ= 31.58981  vector:
\end{verbatim}

\begin{verbatim}
#>  [1]  0.08546276 -0.00288825 -0.09114165 -0.17631489 -0.25552949 -0.32610836 -0.38566627
#>  [8] -0.43219042 -0.46410851 -0.48034187
\end{verbatim}

\begin{verbatim}
#> Minimal eigenvalue = 8.582807e-06  RQ= 8.582807e-06  vector:
\end{verbatim}

\begin{verbatim}
#>  [1] 0.866015080 0.433009398 0.216509345 0.108264428 0.054151958 0.027115582 0.013637056
#>  [8] 0.006977088 0.003805678 0.002537115
\end{verbatim}
\normalsize

First, we will try to find the maximal eigenvalue.

\tiny
\begin{verbatim}
#> Rayleigh Quotient minimization to find maximal eigenvalue
\end{verbatim}

\begin{verbatim}
#> Maximal eigenvalue is calculated as  31.58981  in  1.213167  ms
\end{verbatim}

\begin{verbatim}
#>  [1]  0.08546570 -0.00288988 -0.09114540 -0.17631983 -0.25553608 -0.32611385 -0.38566910
#>  [8] -0.43218968 -0.46410354 -0.48033478
\end{verbatim}

\begin{verbatim}
#> Eigenvalue difference: 5.408658e-09
\end{verbatim}
\normalsize

Then we try to find the minimal eigensolution.

\tiny
\begin{verbatim}
#> Minimal eigenvalue is calculated as  8.582811e-06  in  3.088106  ms
\end{verbatim}

\begin{verbatim}
#>  [1] 0.866015223 0.433008946 0.216509339 0.108265019 0.054151846 0.027116266 0.013637118
#>  [8] 0.006976469 0.003805918 0.002537010
\end{verbatim}

\begin{verbatim}
#> Eigenvalue difference: 4.546869e-12
\end{verbatim}
\normalsize

If we ignore the constraint, and simply perform the optimization, we can sometimes
get satisfactory solutions, though comparisons require that we normalize
the parameters post-optimization. We can check if the scale of the eigenvectors
is becoming large by computing the norm of the final parameter vector. In
tests on the Moler matrix up to dimension 100, none grew to a worrying size.

We note that calculating all the eigensolutions via \texttt{eigen()} takes less time than
a single solution via optimization. For comparison, we also ran a specialized
Geradin routine as implemented in R by
one of us (JN). This gave equivalent answers, albeit more efficiently. For those
interested, the Geradin routine is available as referenced in Nash (\cite{RQtimes12}{2012}).

\section{The extended Rosenbrock function on the unit ball}\label{the-extended-rosenbrock-function-on-the-unit-ball}

The \texttt{adagio} package (Borchers (\cite{p-adagio}{2022})) gives an extended version of the
Rosenbrock banana-shaped valley problem. This becomes a sumscale problem
if we constrain the parameters to be on the unit ball, that is, where
the sum of squares of the parameters is 1.

As with other sumscale problems, the constraint means that while the
Rosenbrock function may be evaluated with \texttt{npar} cartesian parameters,
there are actually only \texttt{npar\ -\ 1} intrinsic parameters. Using this
context, spherical coordinates use a radius \texttt{r=1} and a set of
\texttt{npar\ -\ 1} angles. The Wikipedia article
{https://en.wikipedia.org/wiki/N-sphere} explains these quite well.

However, let us first set up the function in Cartesian coordinates.

\tiny
\begin{verbatim}
# provide a start
x0 <- 1:6/10
library(optimx)
library(BB)
####################################
# Minimizing a function on the unit ball:  Min f(x), s.t. ||x|| = 1
#
rosbkext.f <- function(x){
    p <- x
    n <- length(p)
    sum (100*(p[1:(n-1)]^2 - p[2:n])^2 + (p[1:(n-1)] - 1)^2)
}
\end{verbatim}
\normalsize

We could explicitly constrain the function. The solvers \texttt{alabama::auglag()} and
\texttt{nloptr::slsqp} both allow a vector-valued function \texttt{heq()} such that all elements
of the result of its computation should be 0 on the constraint. That is easy to
specify as follows.

\tiny
\begin{verbatim}
library(alabama)
\end{verbatim}

\begin{verbatim}
#> Loading required package: numDeriv
\end{verbatim}

\begin{verbatim}
library(nloptr)
\end{verbatim}

\begin{verbatim}
#> 
#> Attaching package: 'nloptr'
\end{verbatim}

\begin{verbatim}
#> The following object is masked from 'package:alabama':
#> 
#>     auglag
\end{verbatim}

\begin{verbatim}
heq <- function(x){
    h <- 1 - sum(x*x)
}
\end{verbatim}

\begin{verbatim}
#> Call: ans2 <- alabama::auglag(x0, rosbkext.f, heq=heq,
#>  control.outer = list(trace=FALSE, kktchk=FALSE))
\end{verbatim}

\begin{verbatim}
#> alabama time = 0.014
\end{verbatim}

\begin{verbatim}
#> Result  ans2 (  ->  ) calc. min. = 2.467509  at 
#> 0.748     0.5649892     0.3273791     0.1165197     0.02300184     
#> After  502  fn evals, and  60  gr evals and  NA  hessian evals
#> Termination code is  0 : 
#> Gradient:[1] 0.008002933 0.013716940 0.011464522 0.004793158 0.006851009
#> 
#> -------------------------------------------------
\end{verbatim}

\begin{verbatim}
#> Call: ans3 <- nloptr::slsqp(x0, rosbkext.f, heq=heq)
\end{verbatim}

\begin{verbatim}
#> slsqp time = 0.006
\end{verbatim}

\begin{verbatim}
#> Result  ans3 (  ->  ) calc. min. = 2.467509  at 
#> 0.7480008     0.5649858     0.3273749     0.1165386     0.0230208     
#> After  NA  fn evals, and  45  gr evals and  NA  hessian evals
#> Termination code is  4 : NLOPT_XTOL_REACHED:  
#> 
#> -------------------------------------------------
\end{verbatim}
\normalsize

From the counts, which are incomplete, it seems that \texttt{slsqp} is doing less work
than \texttt{alabama}, but timings are both very short.

We can also project onto the constraint. Let us try with both a scaling of the
parameters and also a projection that uses the sign of the first element.

\tiny
\begin{verbatim}
ProjSphere <- function(x){
    x/sqrt(sum(x*x))
}
\end{verbatim}

\begin{verbatim}
ProjSpheresgn <- function(x){
    sign(x[1])*x/sqrt(sum(x*x))
}
\end{verbatim}

Now we apply \texttt{BB::spg()} to the problem.

\begin{verbatim}
#> ans4 <- spg(x0, rosbkext.f, project=ProjSphere)
\end{verbatim}

\begin{verbatim}
#> iter:  0  f-value:  46.34  pgrad:  1.03722
\end{verbatim}

\begin{verbatim}
#> Warning in spg(x0, rosbkext.f, project = ProjSphere): Unsuccessful convergence.
\end{verbatim}

\begin{verbatim}
#> spg ProjSphere time = 0.147
\end{verbatim}

\begin{verbatim}
#> Result  ans4 (  ->  ) calc. min. = 46.34  at 
#> 0.1048285     0.209657     0.3144855     0.4193139     0.5241424     
#> After  1  fn evals, and  1  gr evals and  NA  hessian evals
#> Termination code is  2 : Maximum function evals exceeded 
#> Gradient:[1] 1.948701       NA       NA       NA       NA
#> 
#> -------------------------------------------------
\end{verbatim}

\begin{verbatim}
#> ans4a <- spg(x0, rosbkext.f, project=ProjSpheresgn)
\end{verbatim}

\begin{verbatim}
#> iter:  0  f-value:  46.34  pgrad:  0.3034194 
#> iter:  10  f-value:  2.849711  pgrad:  1.160761 
#> iter:  20  f-value:  2.467509  pgrad:  0.000501241
\end{verbatim}

\begin{verbatim}
#> spg ProjSpheresgn time = 0.002
\end{verbatim}

\begin{verbatim}
#> Result  ans4a (  ->  ) calc. min. = 2.467509  at 
#> 0.7480006     0.5649858     0.3273752     0.1165392     0.02302091     
#> After  26  fn evals, and  25  gr evals and  NA  hessian evals
#> Termination code is  0 : Successful convergence 
#> Gradient:[1] 5.611599e-06           NA           NA           NA           NA
#> 
#> -------------------------------------------------
\end{verbatim}
\normalsize

The unsigned projection for \texttt{BB::spg()} has failed, while the signed one has
worked reasonably well.

We can also use the angles of spherical coordinates.
If our parameters are the angles (with \texttt{r=1} assumed) of spherical coordinates,
we need to transform them to cartesian coordinates (of which there will be one
more) in order to compute the extended Rosenbrock function on the surface of
the unit ball.

\tiny
\begin{verbatim}
tran <- function(ax){ # NOT right 
# transforms x into z such that ||z|| = 1
    n1 <- length(ax)
    n <- n1 + 1
    z <- rep(NA, n) # this has p+1 elements
    z[1] <- cos(ax[1])
    z[n] <- prod(sin(ax[-n]))
    # should we bound ax with the principal values of inverse trig functions?
    if (n > 2) z[2:(n-1)] <- cumprod(sin(ax[1:(n-2)]))*cos(ax[2:(n-1)])
    return(z)
}
\end{verbatim}
\normalsize

Once we have the \texttt{tran} function, we can use it to convert the angles to the
cartesian parameters and compute the objective
function on the unit ball.

\tiny 
\begin{verbatim}
rosbkang.f <- function(ax){
    n <- length(ax) + 1
    z <- tran(ax) # p are the cartesian coordinates
    ff <- rosbkext.f(z)
}
\end{verbatim}
\normalsize

We need a set of starting angles. Let us try all 1's for \texttt{n=5}.

\tiny 
\begin{verbatim}
p0 <- rep(1,5) # problem in 6 cartesian coordinates
xp0 <- tran(p0)
print(xp0)
\end{verbatim}

\begin{verbatim}
#> [1] 0.5403023 0.4546487 0.3825737 0.3219247 0.2708903 0.4218866
\end{verbatim}

\begin{verbatim}
cat("starting function value at xp0=",rosbkext.f(xp0),"\n")
\end{verbatim}

\begin{verbatim}
#> starting function value at xp0= 25.64734
\end{verbatim}
\normalsize

Let us use this function with some unconstrained solvers to see if we can find
a minimum of the constrained extended Rosenbrock function on the unit ball.

\tiny 
\begin{verbatim}
#> opm(p0, fn=rosbkang.f, gr='grcentral', method=sevmeth, control=list(trace=0))
\end{verbatim}

\begin{verbatim}
#> Warning in BB::spg(par = spar, fn = efn, gr = egr, lower = slower, upper = supper, :
#> Unsuccessful convergence.
\end{verbatim}

\begin{verbatim}
#> Results in spherical coordinates:
\end{verbatim}

\begin{verbatim}
#>                  p5 s5    value fevals gevals hevals conv kkt1 kkt2 xtime
#> nvm      0.02269179    2.467509     45     23      0    0 TRUE TRUE 0.006
#> tnewt    0.02269066    2.467509    127    126      0    0 TRUE TRUE 0.022
#> ncg      0.02269450    2.467509    203     81      0    0 TRUE TRUE 0.024
#> BFGS     0.02268234    2.467509     50     20      0    0 TRUE TRUE 0.004
#> L-BFGS-B 0.02271958    2.467509     42     42      0    0 TRUE TRUE 0.007
#> spg      0.05222296    2.467556   1216   1003      0    1 TRUE TRUE 0.205
\end{verbatim}

\begin{verbatim}
#> Results in cartesian parameters:
\end{verbatim}

\begin{verbatim}
#> ncg  fval= 2.467509  at [1] 0.7480006 0.5649858 0.3273752 0.1165392 0.0230209 0.0005225
#> nvm  fval= 2.467509  at [1] 0.7480006 0.5649858 0.3273752 0.1165392 0.0230209 0.0005225
#> BFGS  fval= 2.467509  at [1] 0.7480006 0.5649859 0.3273751 0.1165390 0.0230209 0.0005223
#> L-BFGS-B  fval= 2.467509  at [1] 0.7480013 0.5649860 0.3273743 0.1165358 0.0230240 0.0005232
#> tnewt  fval= 2.467509  at [1] 0.7480006 0.5649858 0.3273752 0.1165392 0.0230209 0.0005224
#> spg  fval= 2.467556  at [1] 0.747996 0.564985 0.327385 0.116540 0.023017 0.001203
\end{verbatim}
\normalsize

Here \texttt{BB::spg()} is reported as not converged, and we see the objective function
value is a little higher for this solver.

One of the issues with using angles as the parameters for optimization is that
they will give equivalent function values if any is altered by a shift of \texttt{2*pi}.
Therefore, we could bound the angles to be between \(-pi\) and \(pi\).

\tiny 
\begin{verbatim}
#> 
#>  Add bounds [-pi, pi] for spherical angles
\end{verbatim}

\begin{verbatim}
#> Non-bounds methods requested:[1] "BFGS"
\end{verbatim}

\begin{verbatim}
#> Warning in opm(p0, fn = rosbkang.f, gr = "grcentral", lower = -pi, upper = pi, : A method
#> requested does not handle bounds
\end{verbatim}

\begin{verbatim}
#> Warning in BB::spg(par = spar, fn = efn, gr = egr, lower = slower, upper = supper, :
#> Unsuccessful convergence.
\end{verbatim}

\begin{verbatim}
#>                  p5 s5    value fevals gevals hevals conv kkt1 kkt2 xtime
#> nvm      0.02269188    2.467509     48     24      0    0 TRUE TRUE 0.006
#> tnewt    0.02269177    2.467509    144    143      0    0 TRUE TRUE 0.026
#> ncg      0.02269548    2.467509    173     70      0    0 TRUE TRUE 0.015
#> L-BFGS-B 0.02270156    2.467509     51     51      0    0 TRUE TRUE 0.010
#> spg      0.03617910    2.467519   1240   1003      0    1 TRUE TRUE 0.243
\end{verbatim}

\begin{verbatim}
#> Results in cartesian parameters:
\end{verbatim}

\begin{verbatim}
#> ncg  fval= 2.467509  at [1] 0.7480006 0.5649858 0.3273752 0.1165392 0.0230209 0.0005226
#> nvm  fval= 2.467509  at [1] 0.7480006 0.5649858 0.3273752 0.1165392 0.0230209 0.0005225
#> L-BFGS-B  fval= 2.467509  at [1] 0.7480010 0.5649853 0.3273755 0.1165380 0.0230207 0.0005227
#> tnewt  fval= 2.467509  at [1] 0.7480006 0.5649858 0.3273752 0.1165391 0.0230209 0.0005225
#> spg  fval= 2.467519  at [1] 0.7479993 0.5649861 0.3273764 0.1165400 0.0230243 0.0008334
\end{verbatim}
\normalsize

We no longer have results for \texttt{ucminf} or \texttt{optim::BFGS} which do not handle bounds.
\texttt{spg} has once again failed to converge fully, while the other results more or
less match those without bounds applied. However, it is interesting to note that
the KKT checks pass for all cases, so the returned solution for \texttt{spg()} must be
``close''.
\texttt{spg()} works well with a suitable projection, but for some reason does much
less well on the unconstrained minimization of the transformed objective \texttt{rosbkang.f()}.

\section{The R-help example}\label{the-r-help-example}

As a final example, let us use our present techniques to solve the
problem posed by Lanh Tran on R-help. We will use
only a method that scales the parameters directly inside the objective function and
not bother with gradients for this small problem.

\tiny 
\begin{verbatim}
ssums<-function(x){
  n<-length(x)
  tt<-sum(x)
  ss<-1:n
  xx<-(x/tt)*(x/tt)
  sum(ss*xx)
}
\end{verbatim}
\normalsize

We applied ``MOST'' methods available in package \texttt{optimx} solvers. Once the results
were scaled equivalently, all the solutions were essentially the same, and we will
not display them here.

\tiny 
\begin{verbatim}
cat("Try scaled sum\n")
require(optimx)
st<-runif(3)
aos<-opm(st, ssums, gr="grcentral", method="MOST")
# rescale the parameters
nsol<-dim(aos)[1]
for (i in 1:nsol){ 
  tpar<-aos[i,1:3] 
  ntpar<-sum(tpar)
  tpar<-tpar/ntpar
#  cat("Method ",aos[i, "meth"]," gives fval =", ssums(tpar))
  aos[i, 1:3]<-tpar 
}
summary(aos,order=value)
\end{verbatim}
\normalsize

We can also use an unscaled function and a projection with \texttt{spg()}.

\tiny 
\begin{verbatim}
ssum<-function(x){
  n<-length(x)
  ss<-1:n
  xx<-x*x
  sum(ss*xx)
}
proj.simplex <- function(y) {
# project an n-dim vector y to the simplex Dn
# Dn = { x : x n-dim, 1 >= x >= 0, sum(x) = 1}
# Ravi Varadhan, Johns Hopkins University
# August 8, 2012
n <- length(y)
sy <- sort(y, decreasing=TRUE)
csy <- cumsum(sy)
rho <- max(which(sy > (csy - 1)/(1:n)))
theta <- (csy[rho] - 1) / rho
return(pmax(0, y - theta))
}
set.seed(1234)
st<-runif(3)
as<-spg(st, ssum, project=proj.simplex)
\end{verbatim}

\begin{verbatim}
#> iter:  0  f-value:  1.901089  pgrad:  1
\end{verbatim}

\begin{verbatim}
cat("Using project.simplex with spg: fmin=",as$value," at \n")
\end{verbatim}

\begin{verbatim}
#> Using project.simplex with spg: fmin= 0.5454545  at
\end{verbatim}

\begin{verbatim}
print(as$par)
\end{verbatim}

\begin{verbatim}
#> [1] 0.5454579 0.2727254 0.1818168
\end{verbatim}
\normalsize

Apart from the parameter rescaling, this is an entirely ``doable'' problem.
Note that we can also solve the problem as a Quadratic Program using
the \texttt{quadprog} package. An infelicity in this approach is that the
\texttt{solution} parameters are fine, but the minimium is reported as the
\texttt{Lagrangian} value.

\tiny 
\begin{verbatim}
library(quadprog)
Dmat<-diag(c(1,2,3))
Amat<-matrix(c(1, 1, 1), ncol=1)
bvec<-c(1)
meq=1
dvec<-c(0, 0, 0)
ans<-solve.QP(Dmat, dvec, Amat, bvec, meq=0, factorized=FALSE)
ans
\end{verbatim}

\begin{verbatim}
#> $solution
#> [1] 0.5454545 0.2727273 0.1818182
#> 
#> $value
#> [1] 0.2727273
#> 
#> $unconstrained.solution
#> [1] 0 0 0
#> 
#> $iterations
#> [1] 2 0
#> 
#> $Lagrangian
#> [1] 0.5454545
#> 
#> $iact
#> [1] 1
\end{verbatim}
\normalsize

\section{Conclusion}\label{conclusion}

Sumscale problems can present difficulties for optimization (or function minimization)
codes. These difficulties are by no means insurmountable, but they do require some
attention.

While specialized approaches are ``best'' for speed and correctness, finding the
appropriate software may be time consuming, so a general user
is more likely to be able to find solutions quickly from the simpler approach of
embedding the scaling in the
objective function of a conventional optimization solver. In some cases an
unscaled objective where we rescale
the parameters before reporting them also succeeds, but we do not recommend such
approaches because of the burden of checking the proposed solutions. We also
note that the use of a projected gradient via \texttt{spg()} from package \texttt{BB} works
very well, but the projection needs to be set up carefully, as with the use of
the sign of the first element in dealing with the Rayleigh Quotient and the
extended Rosenbrock function on the unit ball.

\section{Acknowledgements}\label{acknowledgements}

We are grateful to Gabor Grothendieck and Hans Werner Borchers for suggestions on
ways to approach this problem.

\section{References}\label{references}

\bibliographystyle{siam}
\bibliography{sumscale.bib}

\section*{Addresses}

\noindent
John C. Nash\\
University of Ottawa (retired)\\%
Telfer School of Management\\ Ottawa, Ontario\\
{ORCiD: {https://orcid.org/0000-0002-2762-8039}}\\%
{mailto:profjcnash@gmail.com}%

\vspace{5mm}

\noindent
Ravi Varadhan\\
Johns Hopkins University\\%
Johns Hopkins University Medical School,\\ Baltimore, Maryland\\
{mailto:ravi.varadhan@jhu.edu}%

\end{document}